\PassOptionsToPackage{unicode}{hyperref}
\PassOptionsToPackage{hyphens}{url}
\documentclass[
]{article}
\usepackage{xcolor}
\usepackage[margin=1in]{geometry}
\usepackage{amsmath,amssymb}
\usepackage{iftex}
\ifPDFTeX
  \usepackage[T1]{fontenc}
  \usepackage[utf8]{inputenc}
  \usepackage{textcomp} 
\else 
  \usepackage{unicode-math} 
  \defaultfontfeatures{Scale=MatchLowercase}
  \defaultfontfeatures[\rmfamily]{Ligatures=TeX,Scale=1}
\fi
\usepackage{lmodern}
\ifPDFTeX\else
\fi
\IfFileExists{upquote.sty}{\usepackage{upquote}}{}
\IfFileExists{microtype.sty}{
  \usepackage[]{microtype}
  \UseMicrotypeSet[protrusion]{basicmath} 
}{}
\makeatletter
\@ifundefined{KOMAClassName}{
  \IfFileExists{parskip.sty}{%
    \usepackage{parskip}
  }{
    \setlength{\parindent}{0pt}
    \setlength{\parskip}{6pt plus 2pt minus 1pt}}
}{
  \KOMAoptions{parskip=half}}
\makeatother
\usepackage{bookmark}
\IfFileExists{xurl.sty}{\usepackage{xurl}}{} 
\hypersetup{
  pdftitle={Connecting Riemannian Geometry and Statistical Inference for Correlation Matrices},
  hidelinks,
  pdfcreator={LaTeX via pandoc}}

\title{Connecting Riemannian Geometry and Statistical Inference for
Correlation Matrices}
\author{Argyn Kuketayev}
\date{\today}

\begin{document}
\maketitle

\subsection{Abstract}\label{abstract}

The quotient-affine metric gives an intrinsic Riemannian geometry to
full-rank correlation matrices, but its geodesic distance has no closed
form and we are not aware of an analytic asymptotic null distribution
for it. We connect this geometry, introduced in 2019, with Jennrich's
1970 asymptotic test for equality of correlation matrices. The quadratic
form underlying Jennrich's statistic is exactly one half of the
quotient-affine metric tensor. The identity arises because eliminating
marginal standard deviations from Gaussian Fisher information performs
the same projection as quotienting out diagonal rescalings. Jennrich's
statistic therefore evaluates the local quotient-affine quadratic form
directly. Moreover, for two independent Gaussian samples with a common
population correlation matrix, the squared geodesic distance, scaled by
effective sample size, converges in distribution to \(4\chi^2_d\), where
\(d = p(p-1)/2\). For \(p=2\), the result reduces to the two-sample
Fisher \(z\) test.

\textbf{Keywords:} correlation matrix; quotient-affine metric; profile
information; Fisher \(z\); asymptotic chi-squared test

\begin{center}\rule{0.5\linewidth}{0.5pt}\end{center}

\subsection{1. Introduction}\label{introduction}

A correlation matrix is a covariance matrix with the marginal scales
discarded. We connect two developments separated by almost fifty years:
Jennrich's 1970 asymptotic test for correlation matrices and the
quotient-affine Riemannian geometry introduced in 2019.

David and Gu (2019) identified the manifold \(\operatorname{Cor}_+(p)\)
of full-rank correlation matrices with the quotient \[
\operatorname{Cor}_+(p) = \operatorname{SPD}(p)/\operatorname{Diag}^+(p),
\] where positive diagonal matrices act on covariance matrices by
congruence, \(\Sigma \mapsto D\Sigma D\). Each fibre contains covariance
matrices with the same correlation structure, so motion along a fibre
changes only marginal scales. Equipping \(\operatorname{SPD}(p)\) with
the affine-invariant metric
\(G_\Sigma(A,B) = \operatorname{tr}(\Sigma^{-1}A\Sigma^{-1}B)\) induces
a quotient metric \(g^Q\) on \(\operatorname{Cor}_+(p)\). Thanwerdas and
Pennec (2021) later derived this metric in closed form.

The associated distance \(d_Q\) is intrinsic, but its computation is
less direct. The Riemannian logarithm has no closed-form expression, so
evaluating \(d_Q(C_1,C_2)\) requires minimizing the affine-invariant
distance over the fibre above \(C_2\); uniqueness of the minimizer is
not known. We are also not aware of an analytic asymptotic null
distribution for this distance.

Jennrich (1970), working from a different direction, derived an
asymptotic chi-squared test for equality of two correlation matrices
under Gaussian sampling. He started from the Gaussian information matrix
in the parameters \((\mu,\sigma,\rho)\), treated the standard deviations
as nuisance parameters, and obtained the information for the
correlations as a Schur complement. His statistic is a quadratic form in
the difference of the two sample correlation matrices. Its first term is
the usual covariance-matrix quadratic form; a second term, described by
Jennrich as ``a correction employed when testing correlation matrices,''
subtracts the contribution of marginal rescaling.

Written in common notation, Jennrich's quadratic form and the
Thanwerdas--Pennec metric are the same expression up to a factor of two
(Proposition 1). The correction term is the squared component associated
with changes of marginal scale, exactly the part removed by the quotient
construction. For the two-sample testing problem, Jennrich's statistic
therefore evaluates the relevant local quotient-affine quadratic form
directly, without computing the geodesic distance. Conversely, the
geodesic distance inherits Jennrich's null law: under Gaussian sampling
with a common population correlation matrix, \[
n_{\mathrm{eff}}\,d_Q^2(R_1,R_2) \xrightarrow{d} 4\chi^2_d,
\qquad d = \frac{p(p-1)}{2},
\] for every \(C \in \operatorname{Cor}_+(p)\).

The passage from the metric tensor to the geodesic distance is
inherently local. Under the null, two sample correlation matrices
approach the same population matrix at the \(n^{-1/2}\) scale, so only
the second-order behaviour of the squared geodesic distance contributes
asymptotically; the global nonlinear distance is replaced, to leading
order, by its tangent-space quadratic form. Proposition 1 identifies
that quadratic form with Jennrich's statistic, allowing his asymptotic
theory to be transferred to the quotient-affine distance.

The contribution is the connection between these developments.
Jennrich's correction term acquires a geometric interpretation as the
component removed when marginal scales are quotiented out, and his
asymptotic sampling theory carries over to the quotient-affine distance.
Jennrich's form is by construction the inverse asymptotic covariance of
the sample correlations, obtained as a Schur complement of the Gaussian
information; Neudecker and Satorra (1996) verified that his expression
agrees algebraically with an independently derived expression for that
variance matrix. The identification here is of a different kind: the
same form is a Riemannian metric tensor. The broader relation between
profile information, Schur complements, and horizontal projection is
well established; see Fewster and Jupp (2013, Remark 1) and
Barndorff-Nielsen and Jupp (1988, §3).

\begin{center}\rule{0.5\linewidth}{0.5pt}\end{center}

\subsection{2. The bivariate case}\label{the-bivariate-case}

The result is visible in the smallest nontrivial case. For \(p=2\) write
\(C(\rho)\) for the correlation matrix with off-diagonal
\(\rho \in (-1,1)\). Thanwerdas and Pennec (2021, Theorem 6) give \[
d_Q\bigl(C(\rho_1),C(\rho_2)\bigr) = \sqrt{2}\,\bigl|z(\rho_2) - z(\rho_1)\bigr|,
\qquad z = \operatorname{arctanh},
\] so the intrinsic distance between two bivariate correlation matrices
is exactly \(\sqrt2\) times the difference of their Fisher \(z\) scores.
With \(z(r_1) - z(r_2) \approx N(0,\,n_1^{-1}+n_2^{-1})\) under a common
\(\rho\), this gives
\(n_{\mathrm{eff}}\,d_Q^2 \xrightarrow{d} 4\chi^2_1\), which is
Corollary 1 at \(p = 2\).

The bivariate case is one-dimensional, and therefore cannot distinguish
the quotient-affine metric from other Riemannian metrics on correlation
matrices; several coincide there. It is an orientation, not a special
case of a transformation. The general statement is Proposition 1, and in
higher dimension there is no analogous scalar transformation: Jennrich's
quadratic form is the corresponding local Fisher-information object.

\begin{center}\rule{0.5\linewidth}{0.5pt}\end{center}

\subsection{3. The identity, and why it
holds}\label{the-identity-and-why-it-holds}

The tangent space at \(C \in \operatorname{Cor}_+(p)\) is the space of
symmetric hollow matrices,
\(T_C\operatorname{Cor}_+(p) = \{Y = Y^\top : \operatorname{diag}(Y) = 0\}\).
Write \[
A(C) = C \circ C^{-1}, \qquad T(C) = I + A(C), \qquad q_C(Y) = \operatorname{diag}(C^{-1}Y) \in \mathbb{R}^p,
\] with \(\circ\) the Hadamard product and
\(\operatorname{diag}(\cdot)\) the diagonal as a column vector.

Theorem 3 of Thanwerdas and Pennec (2021) gives the quotient-affine
metric as \[
g^Q_C(Y,Y) = \operatorname{tr}(C^{-1}YC^{-1}Y) - 2\,q_C(Y)^\top T(C)^{-1} q_C(Y),
\tag{1}
\] and Jennrich's equation (3.9), in the same notation, is \[
\chi^2_J(Y,C) = \tfrac12 \operatorname{tr}(C^{-1}YC^{-1}Y) - q_C(Y)^\top T(C)^{-1} q_C(Y).
\tag{2}
\]

\begin{quote}
\textbf{Proposition 1.} For every \(C \in \operatorname{Cor}_+(p)\) and
every \(Y \in T_C\operatorname{Cor}_+(p)\), \[
\chi^2_J(Y,C) = \tfrac12\, g^Q_C(Y,Y).
\]
\end{quote}

\emph{Proof.} Compare (1) and (2) term by term. \(\square\)

The algebra is immediate; the substance lies in recognizing that the two
formulas describe the same local object.

For a centred Gaussian model \(N(0,\Sigma)\), the Fisher information in
covariance directions is \[
I_\Sigma(A,B) = \tfrac12\operatorname{tr}(\Sigma^{-1}A\Sigma^{-1}B) = \tfrac12 G_\Sigma(A,B),
\tag{3}
\] Jennrich's own equation (3.2): the Gaussian Fisher metric on
covariance matrices is one half of the affine-invariant metric.

Write \(\Sigma = \Delta C \Delta\) with \(\Delta\) positive diagonal,
and parameterize the marginal scales by \(\eta_k = \log\sigma_k\). At
\(\Delta = I\), \[
\frac{\partial \Sigma}{\partial \eta_k} = e_k c_k^\top + c_k e_k^\top =: E_k, \qquad c_k = Ce_k,
\tag{4}
\] which is Jennrich's equation (3.6). Changing a marginal standard
deviation moves the covariance matrix without changing its correlation
matrix. These are exactly the directions removed by the quotient (the
vertical directions). Their inner products under the affine-invariant
metric are
\(\langle E_k,E_l\rangle_C = 2\bigl(\delta_{kl} + C_{kl}(C^{-1})_{kl}\bigr)\),
so Jennrich's matrix \(T(C)\) is, up to the factor two, the matrix that
measures these scale directions.

Jennrich removes the same directions from the Fisher information by the
Schur complement \[
I_{\mathrm{eff}} = I_{CC} - I_{C\eta}I_{\eta\eta}^{-1}I_{\eta C},
\tag{5}
\] his equation (3.4), with \(I_{\eta\eta} = T(C)\) by his (3.7). A
Schur complement of a positive-definite form is orthogonal projection
onto the complement of the eliminated directions; Fewster and Jupp
(2013, Remark 1) give this in coordinate-free form and name the result
the \emph{horizontal information}.

Profiling the marginal scales out of the Gaussian information and
quotienting \(\operatorname{SPD}(p)\) by \(\operatorname{Diag}^+(p)\)
remove the same directions under the same ambient metric. Since the
Fisher metric is \(\tfrac12 G\), the efficient Fisher metric on
correlations is \(\tfrac12 g^Q\), which gives the quadratic-form
identity in Proposition 1.

\begin{center}\rule{0.5\linewidth}{0.5pt}\end{center}

\subsection{4. The geodesic distance}\label{the-geodesic-distance}

Let \(R_1, R_2\) be independent sample correlation matrices from
Gaussian populations with common full-rank \(C\), based on samples of
size \(n_1, n_2\). Put \[
a_n = \frac{n_1n_2}{n_1+n_2}, \qquad n_{\mathrm{eff}} = 2a_n, \qquad \bar R = \frac{n_1R_1 + n_2R_2}{n_1+n_2}.
\] Jennrich's statistic is \(J_n = \chi^2_J(X_n,\bar R)\) with
\(X_n = \sqrt{a_n}\,(R_1-R_2)\), and \(J_n \xrightarrow{d} \chi^2_d\),
\(d = p(p-1)/2\). By Proposition 1, \[
J_n = \tfrac{a_n}{2}\, g^Q_{\bar R}(R_1-R_2,\,R_1-R_2).
\tag{6}
\]

Suppose \(p\) is fixed and \(n_1/(n_1+n_2) \to \eta \in (0,1)\). Under
the null \(R_r - C = O_p(n_r^{-1/2})\). Since the squared Riemannian
distance is smooth near the diagonal, a Taylor expansion about \((C,C)\)
in the matrix-coordinate chart gives \[
d_Q^2(R_1,R_2) = g^Q_C(R_1-R_2,\,R_1-R_2) + O_p(a_n^{-3/2}),
\tag{7}
\] while the leading quadratic term is \(O_p(a_n^{-1})\); hence the
remainder is \(o_p(a_n^{-1})\). Replacing \(C\) by the consistent
\(\bar R\) perturbs the quadratic term by \(o_p(a_n^{-1})\) also, so
combining (6) and (7) gives \(n_{\mathrm{eff}}\,d_Q^2 = 4J_n + o_p(1)\).

\begin{quote}
\textbf{Corollary 1.} For fixed \(p\), under Gaussian sampling with
\(P_1 = P_2 = C \in \operatorname{Cor}_+(p)\), \[
n_{\mathrm{eff}}\, d_Q^2(R_1,R_2) \xrightarrow{d} 4\chi^2_d, \qquad d = \tfrac{p(p-1)}{2}.
\]
\end{quote}

The limit does not depend on \(C\), so the quotient-affine geodesic
distance is asymptotically pivotal. For \(n_1 = n_2 = n\),
\(n_{\mathrm{eff}} = n\). If both populations are elliptical with common
kurtosis parameter \(\kappa\) and finite fourth moments, Neudecker's
(1996) proportionality result replaces the limit by
\(4(1+\kappa)\chi^2_d\).

Proposition 1 is an exact identity between quadratic forms, whereas
Corollary 1 is local and asymptotic: as noted in Section 1, the metric
tensor governs the leading behaviour only because the two sample
matrices approach one another at the \(n^{-1/2}\) scale. For finite or
larger separations, the full nonlinear geometry need not be captured by
Jennrich's quadratic approximation, and outside the elliptical class the
parameter-free limit need not persist.

Higher-order departures from the local quadratic approximation reflect
the nonlinear geometry of the quotient-affine manifold; in normal
coordinates, curvature enters the expansion at fourth order. Closed-form
expressions for this curvature are available from Thanwerdas and Pennec
(2021, 2022), but quantifying the resulting approximation error is a
separate problem.

\begin{center}\rule{0.5\linewidth}{0.5pt}\end{center}

\subsubsection{Declaration of generative AI and AI-assisted technologies
in the manuscript preparation
process.}\label{declaration-of-generative-ai-and-ai-assisted-technologies-in-the-manuscript-preparation-process.}

During the preparation of this work, the author used OpenAI ChatGPT to
assist with literature discovery, manuscript organization, drafting, and
language editing. The author reviewed and edited the output as needed
and takes full responsibility for the content of the published article.

\subsection{References}\label{references}

Barndorff-Nielsen, O. E., \& Jupp, P. E. (1988). Differential geometry,
profile likelihood, \(L\)-sufficiency and composite transformation
models. \emph{The Annals of Statistics}, \textbf{16}(3), 1009--1043.

David, P., \& Gu, W. (2019). A Riemannian structure for correlation
matrices. \emph{Operators and Matrices}, \textbf{13}(3), 607--627.

Fewster, R. M., \& Jupp, P. E. (2013). Information on parameters of
interest decreases under transformations. \emph{Journal of Multivariate
Analysis}, \textbf{120}, 34--39.

Jennrich, R. I. (1970). An asymptotic \(\chi^2\) test for the equality
of two correlation matrices. \emph{Journal of the American Statistical
Association}, \textbf{65}(330), 904--912.

Neudecker, H. (1996). The asymptotic variance matrices of the sample
correlation matrix in elliptical and normal situations and their
proportionality. \emph{Linear Algebra and its Applications},
\textbf{237/238}, 127--132.

Neudecker, H., \& Satorra, A. (1996). The algebraic equality of two
asymptotic tests for the hypothesis that a normal distribution has a
specified correlation matrix. \emph{Statistics \& Probability Letters},
\textbf{30}, 99--103.

Thanwerdas, Y., \& Pennec, X. (2021). Geodesics and curvature of the
quotient-affine metrics on full-rank correlation matrices. In
\emph{Geometric Science of Information (GSI 2021)}, LNCS \textbf{12829},
93--102. Springer.

Thanwerdas, Y., \& Pennec, X. (2022). Theoretically and computationally
convenient geometries on full-rank correlation matrices. \emph{SIAM
Journal on Matrix Analysis and Applications}, \textbf{43}(4),
1851--1872.

\end{document}